\documentstyle[11pt]{article}
\newcommand{\bbz}{\mbox{\boldmath $Z$}}

\newcommand{\bbc}{\mbox{\boldmath $C$}}

\newcommand{\qed}{\hfill$\Box$}

\newcommand{\la}{\lambda}

\newcommand{\Par}{\hbox{Par}}

\newtheorem{thm}{Theorem}[section]

\newtheorem{lem}[thm]{Lemma}

\newtheorem{cor}[thm]{Corollary}

\newcommand{\sym}{\mbox{Sym}}

\begin{document}

\title{Hook Interpolations}
\bibliographystyle{acm}
\author{Ron M. Adin%
\thanks{Department of Mathematics and Computer Science, Bar-Ilan University,
Ramat-Gan 52900, Israel. Email: {\tt radin@math.biu.ac.il} }\ $^\S$ 
\and
Avital Frumkin%
\thanks{Sackler School of Mathematical Sciences, Tel-Aviv University,
Ramat-Aviv, Tel-Aviv 69978, Israel.
 Email: {\tt frumkin@math.tau.ac.il}}
\and Yuval Roichman%
\thanks{Department of Mathematics and Computer Science, Bar-Ilan University,
Ramat-Gan 52900, Israel. Email: {\tt yuvalr@math.biu.ac.il} } 
\thanks{Research was supported in part by 
the Israel Science Foundation, founded by the Israel Academy of Sciences 
and Humanities, and by internal research grants from Bar-Ilan University.}} 
\date{submitted July 30, 2001; revised November 14, 2002}

\maketitle

\begin{abstract}

The hook components of  $V^{\otimes n}$ interpolate between the symmetric power 
$\sym^n(V)$ and the exterior power $\wedge^n(V)$.
When $V$ is the vector space of $k\times m$ matrices over $\bbc$,
we decompose the hook components 
into irreducible $GL_k(\bbc)\times GL_m(\bbc)$-modules. In particular,
 classical  theorems 
are proved as boundary cases.
For the algebra of square matrices over $\bbc$,
a bivariate interpolation is presented and studied.
\end{abstract}

\section{\bf Introduction}

The vector space $M_{k,m}$ of $k\times m$ matrices over $\bbc$
carries a (left) $GL_k(\bbc)$-action and a (right)
$GL_m(\bbc)$-action.
A classical Theorem of Ehresmann \cite{Eh} describes the decomposition of
an exterior power of $M_{k,m}$ into irreducible bimodules.
The symmetric analogue was given later (cf. \cite{Ho}).
See Subsection 2.3 below.

In this paper we present a natural interpolation between these theorems,
in terms of hook components of the $n$-th tensor power of $M_{k,m}$. 
Duality and asymptotics of the decomposition of hook components follow.

Similar methods are then applied to the diagonal two-sided
$GL_k(\bbc)$-action on the vector space of $k\times k$ matrices.
Classical theorems of Thrall \cite{Th} and James \cite{Ja} 
(for the symmetric powers of symmetric matrices), 
and of Helgason \cite{He}, Shimura \cite{Sh} and Howe \cite{H} 
(for the symmetric powers of anti-symmetric matrices)
are extended,  and a bivariate interpolation is presented.

Proofs are obtained using the representation theory of the symmetric and 
hyperoctahedral groups, together with Schur-Weyl duality;
no use is made of highest weight theory.

\subsection{Main Results}

Let $M_{k,m}$ be the vector space of $k\times m$ matrices over $\bbc$. 
The tensor power $M_{k,m}^{\otimes n}$ carries a natural $S_n$-action by 
permuting the factors. This action decomposes the tensor power into 
irreducible $S_n$-modules.
Let $M_{k,m}^{\otimes n}(t)$ be the isotypic component of $M_{k,m}^{\otimes n}$
corresponding to 
the irreducible $S_n$-representation indexed by the hook $(n-t,1^t)$,
where $0\le t\le n-1$.  
This component still carries a $GL_k(\bbc)\times GL_m(\bbc)$-action.

\begin{thm}\label{i1}
Let $\la$ and $\mu$ be partitions of $n$, of lengths at most $k$ and $m$,
respectively.
For every $0\le t\le n-1$,
the multiplicity of the  irreducible $GL_k(\bbc)\times GL_m(\bbc)$-module 
 $V^\la_k\otimes V^\mu_m$ in $M_{k,m}^{\otimes n}(t)$ is
$$
{n-1 \choose t} \sum_{i=0}^{t} (-1)^{t-i} \sigma_{\la,\mu}(i) =
{n-1 \choose t} \sum_{i=t+1}^{n} (-1)^{i-t-1} \sigma_{\la,\mu}(i)
$$
where
$$
\sigma_{\la,\mu}(i) := \sum\limits_{\alpha\vdash n-i, \beta\vdash i} 
c^\la_{\alpha\beta} c^\mu_{\alpha\beta'},
$$
$c^\la_{\alpha\beta}$ are Littlewood-Richardson coefficients, 
and $\beta'$ is the partition conjugate to $\beta$.
\end{thm}

See Theorem \ref{hook1} below;
for definitions and notation see Section 2 below.
Theorem \ref{i1} interpolates between two well-known classical theorems  
(Theorems \ref{Eh} and \ref{Howe} below; see the remark following Theorem \ref{hook1}). 

\medskip

The following corollary generalizes the duality between
Theorem \ref{Eh} and Theorem \ref{Howe}.

\begin{cor}\label{i2} 
Let $\mu\subseteq (m^m)$ and $\la$ be partitions of $n$.
For every $0\le t\le n-1$ the multiplicity of 
$V^\la_k\otimes V^\mu_m$ in $M_{k,m}^{\otimes n}(t)$ is equal to the multiplicity of
$V^\la_k\otimes V^{\mu'}_m$ in $M_{k,m}^{\otimes n}(n-1-t)$.
\end{cor}

See Corollary \ref{dual} below.

\medskip

Let $\lambda$ and $\mu$ be partitions of $n$. Define the {\it distance}
$$
d(\la,\mu):={1\over 2}\sum\limits_i |\la_i-\mu_i|.
$$

\begin{thm}\label{i3}
If $V^\la_k\otimes V^\mu_m$ appears as a factor in $M_{k,m}^{\otimes n}(t)$ 
(for some $0\le t\le n-1$) then
$$
d(\la,\mu)< {km}.
$$
\end{thm}

See Theorem \ref{asym} below.
This shows that, for $V^\la_k\otimes V^\mu_m$
to appear in a hook component,
$\la$ and $\mu$ must be very ``close" to each other
(for $k$ and $m$ fixed, $n$ tending to infinity).

\bigskip

Consider now the vector space $M_{k,k}$ of $k\times k$ square matrices over $\bbc$.
Let $M_{k,k}^{\otimes n}(t,j)$ be the component of $M_{k,k}^{\otimes n}(t)$
consisting of tensors with $j$ skew symmetric and $n-j$ symmetric factors.
$M_{k,k}^{\otimes n}(t,j)$ carries a $GL_k(\bbc)$ two-sided diagonal action.
The following theorem describes its decomposition as a $GL_k(\bbc)$-module.

\begin{thm}\label{i4} 
Let $\la$ be a partition of $2n$ of length at most $k$.
For every $0\le t\le n-1$ and $0\le j\le n$,
the multiplicity of $V_k^\la$ in $M_{k,k}^{\otimes n}(t,j)$ is 
$$
{n-1 \choose t} \sum_{i=0}^{t} (-1)^{t-i} \sigma_{\la}(i,j) =
{n-1 \choose t} \sum_{i=t+1}^{n} (-1)^{i-t-1} \sigma_{\la}(i,j)
$$
where
$$
\sigma_{\la}(i,j) := \sum\limits_{|\alpha|+|\beta|+|\gamma|+|\delta|=n,\ 
|\gamma|+|\delta|=i,\ |\beta|+|\delta|=j}
c_{2\cdot \alpha,(2\cdot \beta)',2*\gamma,(2*\delta)'}^\la ,
$$
and the sum is over all partitions $\alpha,\beta,\gamma,\delta$
with total size $n$ such that $\gamma$ and $\delta$ have distinct parts
and total size $i$, and $\beta$ and $\delta$ have total size $j$.
The operations $*$ and $\cdot$  are defined in Subsection 2.1.
Definition of the (extended) Littlewood-Richardson coefficients is given
in Subsection 2.2.
\end{thm}

See Theorem \ref{hook-sq.1} below.
Theorem \ref{i4}, for $t=0$, interpolates between classical results, 
regarding symmetric powers of the spaces of symmetric
and skew symmetric matrices (Theorems \ref{Ja} and \ref{Shi} below). 
Another boundary case, $t=n$, gives an interpolation between
exterior powers of the same matrix spaces.

\begin{cor}\label{i5}
Let $\la\subseteq (k^k)$ be a partition of $2n$. 
For every $0\le t\le n-1$ and $0\le j\le n$,
the multiplicity
of $V^\la_k$ in $M_{k,k}^{\otimes n}(t,j)$ is equal to the multiplicity of
$V^{\la'}_k$ in $M_{k,k}^{\otimes n}(t,n-j)$.
\end{cor}

See Corollary \ref{hook-sq.2} below.

\section{Background and Notation}

\subsection{Partitions}

Let $n$ be a positive integer. A {\it partition} of $n$
is a vector of positive integers
$\lambda=(\lambda_1,\lambda_2,\ldots,\lambda_k)$,
where 
$\lambda_1\ge\lambda_2\ge\ldots\ge\lambda_k$
and  $\lambda_1+\ldots+\lambda_k=n$.
We denote this by $\la\vdash n$.
The {\it size} of a partition $\la\vdash n$, denoted $|\la|$, is $n$,
and its {\it length}, $\ell(\la)$, is the number of parts.
The empty partition $\emptyset$ has size and length zero: $|\emptyset|=\ell(\emptyset)=0$.
The set of all partitions of $n$ with at most $k$ parts is denoted
by $\Par_k(n)$.
   
For a partition  $\lambda=(\lambda_1,\ldots,\lambda_k)$  define the {\it conjugate partition}
$\lambda'=(\lambda'_1,\dots,\lambda'_t)$ by letting $\lambda'_i$ be 
the number of parts of $\lambda$ that have size at least $i$. 

A partition $\lambda=(\lambda_1,\dots,\lambda_k)$ may be viewed as the subset
$$
\{(i,j)\ |\ 1\le i\le k, 1\le j\le \lambda_i\}\subseteq \bbz^2,
$$
the corresponding {\it Young diagram}. Using this interpretation,
we may speak of the intersection $\lambda\cap\mu$ and the set difference
$\la\setminus \mu$ of any two partitions. The set difference is called
a {\it skew shape}; when $\mu\subseteq \lambda$ it is usually denoted
$\lambda/\mu$.

Let $(k^m):=(k,\dots,k)$ ($m$ equal parts). Thus, for example,
$\la\subseteq (k^m)$ means $\la_1\le k$ and $\la'_1\le m$.

We shall also use the Frobenius notation for partitions, defined as follows:
Let $\la$ be a partition of $n$ and set
$d:=\max\limits\{i\ |\ \la_i-i\ge 0\}$ 
(i.e., the length of the main diagonal in the Young diagram of $\la$).
Then the Frobenius notation for $\la$ is
$(\la_1-1,\dots,\la_d-d\ | \ \la'_1-1,\dots,\la'_d-d)$. 

\smallskip

For any partition $\la=(\la_1,\dots,\la_k)$ of $n$ define the following  doubling operation 
$$
2\cdot \la := (2\la_1,\dots,2\la_k)\vdash 2n .
$$ 
If all the parts of $\la$ are distinct, define also
$$
2 * \la:= (\la_1,\dots\la_k\ |\ \la_1-1,\dots\la_k-1)\vdash 2n
$$
in the Frobenius notation.

\subsection{Representations}

For any group $G$ denote the trivial representation by
$1_G$.
In this paper we shall denote
the irreducible $S_n$-modules (Specht modules) by
$S^\la$, and the irreducible $GL_k(\bbc)$-modules (Weyl modules) by $V^\la_k$.

\medskip

The Littlewood-Richardson coefficients describe the decomposition
of tensor products of Weyl
modules. Let $\mu\vdash t$ and $\nu\vdash n-t$. Then
$$
V^\mu_k\otimes V^\nu_k\cong \bigoplus_{\la\vdash n} c^\la_{\mu,\nu} V^\la_k,
$$
for $k\ge \max\{\ell(\la),\ell(\mu),\ell(\nu)\}$ (and the coefficients
$c^\la_{\mu,\nu}$ are then independent of $k$).

By Schur-Weyl duality they are also the coefficients of the outer
product of Specht modules. Namely,
$$
(S^\mu\otimes S^\nu)\uparrow_{S_t\times S_{n-t}}^{S_n} \cong 
\bigoplus_{\la\vdash n} c^\la_{\mu,\nu} S^\la.
$$
The following identity is well-known: 
For all triples of partitions $\lambda, \mu, \nu$
$$
c^\lambda_{\mu,\nu}=c^{\la'}_{\mu',\nu'}. \leqno{(2.a)}
$$
We shall also use the following notation for 
{\em extended Littlewood-Richardson coefficients} :
$$
c^\la_{\alpha,\beta,\gamma,\delta}:=
\sum\limits_{\mu,\nu} c^\la_{\alpha,\mu}c^\mu_{\beta,\nu}c^\nu_{\gamma,\delta} ;
$$
so that
$$
V^\alpha_k\otimes V^\beta_k\otimes V^\gamma_k\otimes V^\delta_k=
\bigoplus_{\la} c^\la_{\alpha,\beta,\gamma,\delta} V^\la_k.
$$
\bigskip

Let $B_n$ be the Weyl group of type $B$ and rank $n$,
also known as the hyperoctahedral group or the group
of signed permutations. 
A {\it bipartition} of $n$ is an ordered pair $(\mu,\nu)$
of partitions of total size $|\mu|+|\nu|=n$. The irreducible 
characters of $B_n$ are indexed by bipartitions of $n$;
denote by $\chi^{\mu,\nu}$ the character indexed by $(\mu,\nu)$.


\medskip

Consider the following natural embeddings of $S_n$ into $B_n$
and of $B_n$ into $S_{2n}$:
$S_{2n}$ is the group of permutations on $\{-n,\dots,-1,1,\dots,n\}$.
$B_n$ is embedded as the subgroup of all $\pi\in S_{2n}$ satisfying
$\pi(-i)=-\pi(i)$ ($1\le i\le n$).
$S_n$ is embedded as the subgroup of all $\pi\in B_n$ satisfying also
$\pi(i)>0$ ($1\le i\le n$).

The following lemmas, used in Section 5, describe certain induced characters 
via the above embeddings.  
Lemma \ref{d.0} is an immediate consequence of 
\cite[Ch. I \S 7 Ex 4, Ch. I \S 8 Ex 5-6, and Ch. VII (2.4)]{Md}. See also \cite{St}.

\begin{lem}\label{d.0} 
$$
1_{B_n}\uparrow^{S_{2n}}_{B_n}=
\chi^{(n),\emptyset}\uparrow^{S_{2n}}_{B_n}=
\sum\limits_{\la\vdash n}\chi^{2\cdot \la} ; \leqno(a)
$$
$$
\chi^{\emptyset,(n)}\uparrow^{S_{2n}}_{B_n}=
\sum\limits_{\la\vdash n}\chi^{(2\cdot \la)'} ; \leqno(b)
$$
$$
\chi^{(1^n),\emptyset}\uparrow_{B_n}^{S_{2n}}=
\sum\limits_{\la\vdash n}\chi^{2* \la} ;
\leqno(c)
$$
$$
\chi^{\emptyset,(1^n)}\uparrow_{B_n}^{S_{2n}}=
\sum\limits_{\la\vdash n}\chi^{(2* \la)'} ;
\leqno(d)
$$
where the last two sums are over partitions with distinct parts.
\end{lem}

\begin{lem}\label{d.2}
$$
\chi^{(n)}\uparrow^{B_n}_{S_n}=\sum\limits_{i=0}^n \chi^{(i),(n-i)}. \leqno(a)
$$
$$
\chi^{(1^n)}\uparrow^{B_n}_{S_n}=\sum\limits_{i=0}^n 
\chi^{(1^i),(1^{n-i})}. \leqno(b)
$$
\end{lem}

For a proof, see Section 6.1.

The following lemma is a special case of the Littlewood-Richardson rule
for $B_n$, cf. \cite[Lemma 7.1]{Ste2}.

\begin{lem}\label{d.3}
$$
\chi^{(i),(n-i)}=(\chi^{(i),\emptyset}\otimes
\chi^{\emptyset,(n-i)})
\uparrow_{B_i\times B_{n-i}}^{B_n}.
$$
\end{lem}

\subsection{Symmetric and Exterior Powers of Matrix Spaces}

In this subsection we cite well-known classical theorems, 
concerning  the decomposition into irreducibles
of symmetric and exterior powers
of matrix spaces,
which are to be generalized in this paper.

\medskip

 Let $M_{k,m}$ be the vector space of  $k\times m$ matrices over $\bbc$.
Then $M_{k,m}$ carries a (left) $GL_k(\bbc)$-action and a (right)
$GL_m(\bbc)$-action.
A classical Theorem of Ehresmann \cite{Eh} (see also \cite{Ko}) describes the decomposition of
an exterior power of $M_{k,m}$ into irreducible $GL_k(\bbc)\times GL_m(\bbc)$-modules.

\begin{thm}\label{Eh}
The $n$-th exterior power of $M_{k,m}$ is isomorphic, as
a $GL_k(\bbc)\times GL_m(\bbc)$-module, to
$$
\wedge^n (M_{k,m}) \cong 
\bigoplus_{\la\vdash n \hbox{ and } \la\subseteq (m^k)} V^\la_k \otimes V^{\la'}_m, 
$$
where $\la'$ is the partition conjugate to $\la$.
\end{thm}

\medskip

The following three results on symmetric powers
were proved several times independently; 
these results
may be found in \cite{Ho} and \cite{GW}.

The symmetric analogue of Theorem \ref{Eh}
was studied, for example, in 
 \cite[(11.1.1)]{Ho}  and \cite[Theorem 5.2.7]{GW}.

\begin{thm}\label{Howe}
The $n$-th symmetric power of $M_{k,m}$ is isomorphic, as
a $GL_k(\bbc)\times GL_m(\bbc)$-module, to
$$
\sym^n (M_{k,m})\cong 
\bigoplus_{\la\vdash n \hbox{ and } \ell(\la)\le \min(k,m)} V^\la_k \otimes V^\la_m. 
$$
\end{thm}

Let $M_{k,k}^+$ be the vector space of symmetric $k\times k$ matrices
over $\bbc$.
This space carries a natural two sided $GL_k(\bbc)$-action. 
The following 
theorem describes the decomposition of its symmetric powers
into irreducible $GL_k(\bbc)$-modules.

\begin{thm}\label{Ja}
 The $n$-th symmetric power of $M_{k,k}^+$ is isomorphic, as
a $GL_k(\bbc)$-module, to
$$
\sym^n (M_{k,k}^+)\cong \bigoplus_{\la\in \hbox{Par}_k(n) } V^{2\cdot \la}_k.
$$
\end{thm}

This theorem was proved by A.T. James \cite{Ja},
but had already appeared in an early work of Thrall \cite{Th}.
See also 
\cite{H}, \cite{Sh}, \cite[(11.2.2)]{Ho}
and \cite[Theorem 5.2.9]{GW}
 for further proofs and references. 

\medskip

Let $M^-_{k,k}$ be the vector space of skew symmetric $k\times k$ matrices
over $\bbc$.
Then

\begin{thm}\label{Shi}
 The $n$-th symmetric power of $M^-_{k,k}$ is isomorphic, as
a $GL_k(\bbc)$-module, to
$$
\sym^n (M^-_{k,k})\cong \bigoplus_{(2\cdot\la)'\in \hbox{Par}_k(2n) } V^{(2\cdot \la)'}_k.
$$
\end{thm}

This theorem is proved in  \cite{He}, \cite{H}, \cite{Sh}.
See also \cite[(11.3.2)]{Ho} and \cite[Theorem 5.2.11]{GW}.

\section{Hook Components of $M_{k,m}^{\otimes n}$}

Consider $M=M_{k,m}=\bbc^{k\times m}$, the vector space of $k\times m$
matrices over $\bbc$. Then $M\cong V\otimes W$, where $V\cong \bbc^k$ and 
$W\cong\bbc^m$. Thus $M$ carries a (left) $GL(V)$-action and a (right)
$GL(W)$-action, which commute.
Its tensor power $M^{\otimes n} \cong V^{\otimes n}\otimes W^{\otimes n}$
thus carries a $GL(V)\times S_n\times S_n\times GL(W)$ linear representation; 
one copy of the symmetric group $S_n$ permutes the factors in $V^{\otimes n}$, 
and the other copy of $S_n$ permutes the factors in $W^{\otimes n}$. 
The actions of all four groups clearly commute.
We are interested in the $GL(V)\times S_n\times GL(W)$-action on 
$M^{\otimes n}$ obtained through the diagonal embedding
$S_n\hookrightarrow S_n\times S_n$, $\pi \mapsto (\pi,\pi)$.

\begin{lem}\label{diagonal}
$$
M^{\otimes n} \cong \bigoplus_{\la\in \Par_k(n)\ ,\ \nu\in\Par(n)\hbox{ , }\mu\in\Par_m(n)} 
\alpha_{\la\mu\nu} V^\la_k\otimes S^{\nu}\otimes V^\mu_m,
$$
where
$$
\alpha_{\la\mu\nu}:=\langle \chi^\la\chi^\mu\chi^\nu, 1_{S_n}\rangle.
$$
\end{lem}

\noindent{\bf Proof.}
By Schur-Weyl duality (the double commutant theorem) \cite[Theorem 9.1.2]{GW} 
$$
V^{\otimes n}\cong \bigoplus_{\la\in \Par_k(n)} V^\la_k\otimes S^\la
$$
as $GL(V)\times S_n$-modules.
Similarly,
$$
W^{\otimes n}\cong \bigoplus_{\la\in \Par_m(n)} V^\la_m\otimes S^\la
$$
as $GL(W)\times S_n$-modules.
Therefore
$$
M^{\otimes n}\cong V^{\otimes n}\otimes W^{\otimes n}\cong
\bigoplus_{\la\in \Par_k(n)\hbox{ and } \mu\in \Par_m(n)} 
V^\la_k\otimes S^\la\otimes S^\mu \otimes V_m^\mu
$$
as $GL(V)\times S_n\times S_n\times GL(W)$-modules.

Using the diagonal embedding $S_n\hookrightarrow S_n\times S_n$,
$$
M^{\otimes n}\cong 
\bigoplus_{\la\in \Par_k(n)\hbox{ and } \mu\in \Par_m(n)} 
V^\la_k\otimes (S^\la\otimes S^\mu)\downarrow^{S_n\times S_n}_{S_n} \otimes V_m^\mu
$$
as $GL(V)\times S_n\times GL(W)$-modules.

Note that the $S_n$-character of
$(S^\la\otimes S^\mu)\downarrow^{S_n\times S_n}_{S_n}$
is the standard inner tensor product (sometimes called Kronecker
product) of the $S_n$-characters $\chi^\la$ and $\chi^\mu$.
Hence, by elementary representation theory
$$
(S^\la\otimes S^\mu)\downarrow^{S_n\times S_n}_{S_n}\cong 
\bigoplus_{\nu\vdash n} \alpha_{\la\mu\nu} S^\nu,
$$
where
$$
\alpha_{\la\mu\nu} = \langle \chi^\la \chi^\mu,\chi^\nu \rangle =
{1\over n!}\sum\limits_{\pi\in S_n} \chi^\la(\pi)\chi^\mu(\pi)\chi^\nu(\pi) =
\langle \chi^\la\chi^\mu\chi^\nu, 1_{S_n}\rangle.
$$
\qed

In particular, Lemma \ref{diagonal} gives Theorem \ref{Eh} and Theorem \ref{Howe}.

\medskip

\begin{cor}\label{e-h}
$$
\sym^n (M)\cong \bigoplus_{\la\vdash n \hbox{ and } \ell(\la)\le \min(k,m)} 
V^\la_k \otimes V^\la_m. \leqno(1)
$$
$$
\wedge^n (M) \cong \bigoplus_{\la\vdash n \hbox{ and } \la\subseteq (m^k)} 
V^\la_k \otimes V^{\la'}_m. \leqno(2)
$$
\end{cor}

\noindent{\bf Proof.}
$\sym^n(M)$ is the isotypic component of $M^{\otimes n}$ corresponding to 
the trivial character $\chi^{(n)}$ of the symmetric group.
Thus, by Lemma \ref{diagonal}
$$
\sym^n(M)\cong \bigoplus_{\la\in \Par_k(n)\hbox{ and }\mu\in\Par_m(n)} 
\alpha_{\la,\mu,(n)}V^\la_k\otimes S^{(n)}\otimes V^\mu_m.
$$
But by the orthonormality of irreducible characters
$$
\alpha_{\la,\mu,(n)}=\langle \chi^\la\chi^\mu,\chi^{(n)}\rangle=
\langle \chi^\la, \chi^\mu\chi^{(n)}\rangle=
\langle \chi^\la,\chi^\mu\rangle=\delta_{\la\mu}.
$$
This proves (1), namely Theorem \ref{Howe}.

The $n$-th exterior power is the isotypic
component of $M^{\otimes n}$ corresponding to
the sign character $\chi^{(1^n)}$ of the symmetric group. 
Recall that for any partition $\mu\vdash n$, $\chi^\mu\chi^{(1^n)}=\chi^{\mu'}$.
Thus
$$
\alpha_{\la,\mu,(1^n)}=\langle \chi^\la\chi^\mu,\chi^{(1^n)}\rangle=
\langle \chi^\la, \chi^\mu\chi^{(1^n)}\rangle=
\langle \chi^\la,\chi^{\mu'}\rangle=\delta_{\la\mu'}.
$$
This proves (2), namely Theorem \ref{Eh}.

\qed

\bigskip

Let $M$ be the vector space of $k\times m$ matrices as before. 
The tensor power $M^{\otimes n}$ carries a natural $S_n$-action by 
permuting the factors. This action decomposes into irreducible $S_n$-representations.
Let $M^{\otimes n}(t)$ be the component of $M^{\otimes n}$, corresponding to 
the irreducible hook representation $(n-t,1^t)$, $0\le t\le n-1$.  
This component carries a $GL_k(\bbc)\times GL_m(\bbc)$-action.

\begin{thm}\label{hook1}
Let $\la\in \Par_k(n)$ and $\mu\in \Par_m(n)$.
For every $0\le t\le n-1$,
the multiplicity of the irreducible $GL_k(\bbc)\times GL_m(\bbc)$-module 
 $V^\la_k\otimes V^\mu_m$ in $M^{\otimes n}(t)$  is
$$
{n-1 \choose t} \sum_{i=0}^{t} (-1)^{t-i} \sigma_{\la,\mu}(i) =
{n-1 \choose t} \sum_{i=t+1}^{n} (-1)^{i-t-1} \sigma_{\la,\mu}(i)
$$
where
$$
\sigma_{\la,\mu}(i) := \sum\limits_{\alpha\vdash n-i, \beta\vdash i} 
c^\la_{\alpha\beta} c^\mu_{\alpha\beta'},
$$
$c^\la_{\alpha\beta}$ are Littlewood-Richardson coefficients, 
and $\beta'$ is the partition conjugate to $\beta$.
\end{thm}

\noindent{\bf Remark}.
Theorem \ref{hook1} may be considered as an interpolation 
between Theorem \ref{Eh} and Theorem \ref{Howe}. 
Taking $t=0$ gives $M^{\otimes n}(0)\cong \sym^n(M)$, 
and $\beta\vdash 0$ means $\beta=\emptyset$. 
Hence $\la=\alpha=\mu$, yielding multiplicity $\delta_{\la\mu}$. 
This gives Theorem \ref{Howe}.

\noindent
Similarly, taking $t=n-1$ gives $M^{\otimes n}(n-1)\cong \wedge^n(M)$, 
and $\alpha\vdash 0$ means $\alpha=\emptyset$. 
Hence $\la=\beta=\mu'$, yielding multiplicity $\delta_{\la\mu'}$. 
This gives Theorem \ref{Eh}.

\medskip

\noindent{\bf Proof.}
By Lemma \ref{diagonal}
$$
M^{\otimes n}(t)\cong \bigoplus_{\la\in \Par_k(n)\hbox{ and }\mu\in\Par_m(n)} 
\alpha_{\la,\mu,(n-t,1^t)}V^\la_k\otimes S^{(n-t,1^t)}\otimes V^\mu_m
$$
is the decomposition of this component into irreducibles.

Denote by $1_t$ and $\varepsilon_t$ the trivial and sign characters, respectively, of $S_t$.
By the combinatorial interpretation of the
Littlewood-Richardson rule (cf. \cite[Theorem 2.8.13]{JK}), 
for every $0\le t \le n$
$$
(1_{n-t}\otimes \varepsilon_t)\uparrow_{S_{n-t}\times S_t}^{S_n}=
\chi^{(n-t,1^t)}+ \chi^{(n-t+1,1^{t-1})}. \leqno{(3.a)}
$$
Hence, by Frobenius reciprocity
$$
\alpha_{\la,\mu,(n-t,1^t)}+\alpha_{\la,\mu,(n-t+1,1^{t-1})}=
\langle \chi^\la \chi^\mu, \chi^{(n-t,1^t)}+\chi^{(n-t+1,1^{t-1})}\rangle=
$$
$$
=\langle \chi^\la \chi^\mu, (1_{n-t}\otimes \varepsilon_t)\uparrow_{S_{n-t}\times S_t}^{S_n} \rangle
=\langle (\chi^\la \chi^\mu)\downarrow_{S_{n-t}\times S_t}^{S_n}, 
1_{n-t}\otimes \varepsilon_t \rangle=
$$
$$
=\langle \chi^\la \downarrow_{S_{n-t}\times S_t}^{S_n}, 
\chi^\mu\downarrow_{S_{n-t}\times S_t}^{S_n}\cdot (1_{n-t}\otimes \varepsilon_t) \rangle.
$$
By the Littlewood-Richardson rule the last expression is equal to
$$
\langle \sum\limits_{\alpha\vdash n-t, \beta\vdash t} c^\la_{\alpha\beta}
\chi^\alpha\otimes \chi^\beta, \sum\limits_{\alpha\vdash n-t, \beta\vdash t} 
c^\mu_{\alpha\beta} \chi^\alpha\otimes \chi^\beta \cdot
(1_{n-t}\otimes \varepsilon_t)\rangle=
$$
$$
=\langle \sum\limits_{\alpha\vdash n-t, \beta\vdash t} c^\la_{\alpha\beta}
\chi^\alpha\otimes \chi^\beta, \sum\limits_{\alpha\vdash n-t, \beta\vdash t} 
c^\mu_{\alpha\beta} \chi^\alpha\otimes \chi^{\beta'} \rangle=
\sum\limits_{\alpha\vdash n-t, \beta\vdash t} c^\la_{\alpha\beta}
c^\mu_{\alpha\beta'},
$$
which was denoted $\sigma_{\la,\mu}(t)$ in the statement of the theorem.
Alternating summation and the well-known fact
$$
\dim S^{(n-t,1^t)} = {n-1 \choose t}
$$
now complete the proof.
\qed

\medskip

The following corollary generalizes the ``duality" of Theorem \ref{Eh}
and Theorem \ref{Howe}.

\begin{cor}\label{dual} 
Let $\la\in \Par_k(n)$, and let
$\mu,\mu'\in \Par_m(n)$ be conjugate partitions.
Then, for every $0\le t\le n-1$, the multiplicity of 
$V^\la_k\otimes V^\mu_m$ in $M^{\otimes n}(t)$ is equal to the multiplicity of
$V^\la_k\otimes V^{\mu'}_m$ in $M^{\otimes n}(n-1-t)$.
\end{cor}

\noindent{\bf Proof.}
By Theorem \ref{hook1}, we need to show that
$$
{n-1 \choose t} \sum_{i=0}^{t} (-1)^{t-i} \sigma_{\la,\mu}(i) =
{n-1 \choose n-1-t} \sum_{j=n-t}^{n} (-1)^{j-n+t} \sigma_{\la,\mu'}(j).
$$
This follows from
$$
\sigma_{\la,\mu}(i) = \sigma_{\la,\mu'}(n-i),
$$
which in turn follows from (2.a).

\qed

\medskip

\noindent{\bf Examples.} Let $\la\in \Par_k(n),\ \mu,\mu'\in\Par_m(n)$.
The multiplicities of $V^\la_k\otimes V^\mu_m$ in $M^{\otimes n}(t)$
for $t=0$ and $t=n-1$ are given by Theorems \ref{Howe} and \ref{Eh}.
Consider two other pairs of $t$-values.
\begin{itemize}
\item $t=1$:
For $\la=\mu$ the multiplicity is $n-1$ times the number of (inner) corners 
in $\la$, minus 1.
For $\la\not=\mu$ it is $n-1$ if $|\la\setminus \mu|=1$, and zero otherwise.

\item $t=n-2$:
For $\la=\mu'$ the multiplicity is $n-1$ times the number of (inner) corners 
in $\la$, minus 1.
For $\la\not=\mu'$ it is $n-1$ if $|\la\setminus \mu'|=1$, and zero otherwise.

\item $t=2$ ($n>2$): 
For $\la=\mu$ the multiplicity is nonzero iff $\la$ has at least 3 inner corners.
For $\la\ne\mu$ it is nonzero iff there is a partition $\alpha$ of $n-2$ such that
$\la/\alpha$ is a horizontal strip and $\mu/ \alpha$ is a vertical strip,
or vice versa.

\item $t=n-3$ ($n>2$):
Analogous to the previous case, with $\mu$ replaced by $\mu'$.

\end{itemize}


\section{Asymptotics}

Let $\lambda$ and $\mu$ be partitions of the same number
$n$. Recalling the definition of the set difference $\la\setminus \mu$
from Subsection 2.1,
define the {\it distance}
$$
d(\la,\mu):={|\la\setminus \mu|} 
\ (\ ={1\over 2}\sum\limits_i |\la_i-\mu_i|\ )\ .
$$

\begin{lem}\label{dist}
If $V^\la_k\otimes V^\mu_m$ appears as a factor in $M^{\otimes n}(t)$
(for some $0\le t\le n-1$)
then \ $d(\la,\mu)\le t$\ and \ $d(\la,\mu')\le n-1-t$.
\end{lem}

\noindent{\bf Proof.} 
By Theorem~\ref{hook1},
if $V^\la_k\otimes V^\mu_m$ appears as a factor in $M^{\otimes n}(t)$ 
then there exists a pair of partitions, $\alpha\vdash n-i$ and $\beta\vdash i$,
with $0\le i\le t$, such that $c^\la_{\alpha\beta} c^\mu_{\alpha\beta'}\not=0$.
$c^\la_{\alpha\beta}\not= 0 \Rightarrow \alpha\subseteq \la$,
and $c^\mu_{\alpha\beta'}\not=0 \Rightarrow \alpha\subseteq \mu$.
Hence 
$$
|\la\setminus\mu| \le |\la\setminus \alpha| = i \le t.
$$
Using the second expression in Theorem~\ref{hook1}, 
if $V^\la_k\otimes V^\mu_m$ appears as a factor in $M^{\otimes n}(t)$ 
then there exists a pair of partitions, $\alpha\vdash n-j$ and $\beta\vdash j$,
with $t+1\le j\le n$, such that $c^\la_{\alpha\beta} c^\mu_{\alpha\beta'}\not=0$.
$c^\la_{\alpha\beta}\not= 0 \Rightarrow \beta\subseteq \la$,
and $c^\mu_{\alpha\beta'}\not=0 \Rightarrow \beta\subseteq \mu'$.
Hence 
$$
|\la\setminus\mu'| \le |\la\setminus\beta| = n-j \le n-t-1.
$$
Altogether, we get the desired claim.

\qed

Let $\psi$ be an $S_n$-character (not necessarily irreducible).
Define the {\it height} of $\psi$ by
$$
\hbox{height} (\psi) := \max\{\ell(\nu)\ |\ \nu\vdash n,\ 
\langle \psi, \chi^\nu\rangle \not= 0\}.
$$
The following result was proved by Regev.

\begin{lem}\label{height} \cite[Theorem 12]{Re}
For any $\lambda,\mu\vdash n$,
$$
\hbox{height}(\chi^\la \chi^\mu)\le \ell(\la)\cdot\ell(\mu).
$$
\end{lem}

\begin{thm}\label{asym}
If $V^\la_k\otimes V^\mu_m$ appears as a factor in $M^{\otimes n}(t)$
(for some $0\le t\le n-1$) then
$$
d(\la,\mu)< {km}.
$$
\end{thm}

\noindent{\bf Proof.} 
$$
d(\la,\mu) \stackrel{(1)}{\le} t \stackrel{(2)}{\le} 
\hbox{height}(\chi^\la\chi^\mu)-1\stackrel{(3)}{\le}
\ell(\la)\cdot \ell(\mu)-1\le km-1.
$$
Inequalities (1), (2) and (3) follow from Lemmas \ref{dist},
\ref{diagonal} (for $\nu=(n-t,1^t)$)
and \ref{height}, respectively.
\qed 

\medskip

Let $\psi$ be an $S_n$-character (not necessarily irreducible).
Define the {\it width} of $\psi$ by
$$
\hbox{width} (\psi) := \max\{\nu_1\ |\ \nu\vdash n,\ 
\langle \psi, \chi^\nu\rangle \not= 0\}.
$$
The following result of Dvir strengthens Lemma \ref{height}.

\begin{lem}\label{Dvir}\cite[Theorem 1.6]{Dv}
For any $\lambda,\mu\vdash n$,
$$
\hbox{width}(\chi^\la \chi^\mu)=|\la\cap\mu|\leqno(1)
$$
and
$$
\hbox{height}(\chi^\la \chi^\mu)=|\la\cap \mu'|.\leqno(2)
$$
\end{lem}

\smallskip

This result gives another way of proving Theorem \ref{asym}.

\medskip

\noindent{\bf Second Proof of Theorem \ref{asym}.}
$$
d(\la,\mu)=|\la\setminus \mu|=n-|\la\cap\mu|
\stackrel{(1)}{\le} t\stackrel{(2)}{\le} \hbox{height}(\chi^\la\chi^\mu)-1\stackrel{(3)}{=}
|\la\cap\mu'|-1\le km-1.
$$
Inequality (1) follows from Lemma \ref{Dvir}(1),
since $n-t\le\hbox{ width}(\chi^\la\chi^\mu)$. 
Inequality (2) follows from Lemma \ref{diagonal}.
Equality (3)
is Lemma \ref{Dvir}(2).

\qed

\medskip

\noindent{\bf Note :} 
For any two partitions $\la,\mu$ of $n$
with $\ell(\la)\le k$ and $\ell(\mu)\le m$,
$V^\la_k\otimes V^\mu_m$ appears as a factor in $M^{\otimes n}$.
Theorem \ref{asym} shows that, 
in order to appear in a hook component,
$\la$ and $\mu$ must be very ``close" to each other
(for $k$ and $m$ fixed, $n$ tending to infinity).

\section{Square Matrices}

Consider now the vector space $M_k=M_{k, k}$ of $k\times k$ matrices over $\bbc$.
This space carries a diagonal
(left and right) $GL_k(\bbc)$-action, defined by
$$
g(m):= g\cdot m\cdot g^t \qquad 
(\forall g\in GL_k(\bbc)\ , \forall m \in M_{k}).
$$

\subsection{Symmetric Powers}

Recall from Section 2.1 the definition of $2\cdot \la$, for a partition $\la$.

\begin{thm}\label{d.6} 
For $\la\in\Par_k(2n)$,
the multiplicity of $V_k^\la$ in $\sym^{n}(M_k)$ is 
$$
\sum\limits_{|\mu|+|\nu|=n} c_{2\cdot \mu,(2\cdot \nu)'}^\la .
$$
\end{thm}

\begin{cor}\label{d.7}
Let $\la\in \Par(2n)$, $\la\subseteq (k^k)$ (i,e., $\la,\la'\in \Par_k(2n)$). 
Then the multiplicities
of $V^\la_k$ and of $V^{\la'}_k$
in $\sym^n(M_k)$ are equal.
\end{cor}

\noindent{\bf Proof.} This is an immediate consequence of Theorem
\ref{d.6}, applying identity (2.a).

\qed

\bigskip

\noindent{\bf Proof of Theorem \ref{d.6}.}
Let $V\cong \bbc^k$. Then $V\otimes V$ carries a diagonal (left)
$GL_k$-action, and
$$
M_k\cong V\otimes V
$$
as $GL_k$-modules. Thus
$$
M_k^{\otimes n} \cong V^{\otimes 2n}
$$
as $GL_k$-modules. Moreover, $V^{\otimes 2n}$ carries an $S_{2n}\times GL_k$-action :
$S_{2n}$ permutes the $2n$ factors in the tensor product, 
and $GL_k$ acts on all of them simultaneously (on the left). 
The $S_{2n}$- and $GL_k$- actions satisfy Schur-Weyl duality 
(the double commutant theorem), so that
$$
V^{\otimes 2n}\cong \bigoplus_{\la\in \Par_k(2n)} V_k^\la \otimes S^\la,
$$
as $GL_k\times S_{2n}$-modules.

\medskip

Now, $\sym^n(M_k)$ is the part of $M_k^{\otimes n}$ which is invariant 
under the action of $S_n \hookrightarrow S_{2n}$, where the embedding
 $S_n\hookrightarrow S_n\times S_n \subseteq S_{2n}$ is diagonal:
$\pi \longmapsto (\pi,\pi)\ $.
It follows that the multiplicity of $V^\la_k$ in $\sym^n(M_k)$
is equal to the multiplicity of the trivial character $1_{S_n}$ in the
restriction $\chi^\la\downarrow^{S_{2n}}_{S_n}$, where $S_n$
is diagonally embedded.

By Frobenius reciprocity,
$$
\langle 1_{S_n}, \chi^\la\downarrow^{S_{2n}}_{S_n}\rangle=
\langle 1_{S_n}\uparrow^{S_{2n}}_{S_n}, \chi^\la\rangle .
$$
We conclude that, for $\la\in\Par_k(2n)$, 
the multiplicity of $V_k^\la$ in $\sym^{n}(M_k)$ is 
$$
\langle 1_{S_n}\uparrow^{S_{2n}}_{S_n}, \chi^\la\rangle .
$$
We shall compute these multiplicities in several steps.

First, we induce in two steps:
$$
1_{S_n}\uparrow^{S_{2n}}_{S_n} = 
(1_{S_n}\uparrow^{B_n}_{S_n})\uparrow_{B_n}^{S_{2n}}.
$$
By Lemmas \ref{d.2}(a) and \ref{d.3},
$$
(\chi^{(n)}\uparrow^{B_n}_{S_n})\uparrow^{S_{2n}}_{B_n}=
\sum\limits_{i=0}^n \chi^{(i),(n-i)}\uparrow^{S_{2n}}_{B_n}=
$$
$$
=\sum\limits_{i=0}^n 
(\chi^{(i),\emptyset}\otimes
\chi^{\emptyset,(n-i)})\uparrow^{B_{n}}_{B_i\times B_{n-i}}
\uparrow^{S_{2n}}_{B_n} = 
$$
$$
=\sum\limits_{i=0}^n 
(\chi^{(i),\emptyset}\otimes
\chi^{\emptyset,(n-i)})\uparrow^{S_{2n}}_{B_i\times B_{n-i}}. 
$$

Again, let us induce in two steps :
$$
(\chi^{(i),\emptyset}\otimes
\chi^{\emptyset,(n-i)})
\uparrow^{S_{2n}}_{B_i\times B_{n-i}}=
(\ (\chi^{(i),\emptyset}\otimes
\chi^{\emptyset,(n-i)})
\uparrow_{B_i\times B_{n-i}}^{S_{2i}\times S_{2n-2i}}) \uparrow_{S_{2i}\times S_{2n-2i}}
^{S_{2n}}=
$$
$$
=(\ \chi^{(i),\emptyset}\uparrow_{B_i}^{S_{2i}}  \otimes \ 
\chi^{\emptyset,(n-i)}
\uparrow_{B_{n-i}}^{S_{2n-2i}}) \uparrow_{S_{2i}\times S_{2n-2i}}
^{S_{2n}}.
$$
By Lemma \ref{d.0} (a)-(b), the right hand side is equal to
$$
(\sum\limits_{\mu\vdash i} \chi^{2\cdot \mu} \otimes 
\sum\limits_{\nu\vdash n-i} \chi^{(2\cdot \nu)'})
\uparrow_{S_{2i}\times S_{2n-2i}}
^{S_{2n}}.
$$
We conclude that 
$$
 1_{S_n}\uparrow^{S_{2n}}_{S_n}=
\sum\limits_{i=0}^n \ \sum\limits_{\mu\vdash i, \nu\vdash n-i}
(\chi^{2\cdot \mu} \otimes  \chi^{(2\cdot \nu)'})
\uparrow_{S_{2i}\times S_{2n-2i}}
^{S_{2n}}.
$$
Applying the
Littlewood-Richardson rule completes the proof.

\qed

\subsection{A Graded Refinement of Symmetric Powers}

The space $M_k^{\otimes n}$ carries not only an $S_n$-action but also
a $B_n$-action, where the signed permutation $(i,-i)$ $(1\le i\le n)$ acts by 
transposing the $i$-th factor in the tensor product of $n$ square matrices.
$M_k=M_k^+\oplus M_k^-$, where $M_k^+$ ($M_k^-$) is the vector space of 
symmetric (skew symmetric) matrices of order $k\times k$. Consequently,
$M_k^{\otimes n}$ is graded by the number of skew symmetric factors.
The component of $M_k^{\otimes n}$ with $i$ skew symmetric factors, 
denoted $M_k^{\otimes n}(i)$, is invariant under the $B_n$-action, 
as well as under the diagonal two-sided $GL_k$-action. 

\begin{lem}\label{r.3} 
If the irreducible $B_n$-character $\chi^{\mu,\nu}$ appears in the
decomposition of the $B_n$-action on $M_k^{\otimes n}(i)$, then
$|\nu|=i$.
\end{lem}

For a proof see Section 6.2.

\medskip

Since the components $M_k^{\otimes n}(i)$ are invariant under
the $S_n$-action, the $S_n$-invariant subspace
$\sym^n(M_k)$ inherits the grading by the number
of skew symmetric factors.
Let $\sym^n_i(M_k)$ denote the component with
$i$ skew symmetric factors.
The following theorem refines Theorem \ref{d.6}.

\begin{thm}\label{r.1}
For $\la\in\Par_k(2n)$,
the multiplicity of $V_k^\la$ in $\sym^n_i(M_k)$ is 
$$
\sum\limits_{\mu\vdash n-i\ ,\ \nu\vdash i} c_{2\cdot \mu,(2\cdot \nu)'}^\la .
$$
\end{thm}

\noindent{\bf Note:}
Theorem \ref{r.1} interpolates between two classical results,
Theorem \ref{Ja}
and Theorem \ref{Shi}. Indeed,
$\sym^n_0(M_k)= \sym^n(M_k^+)$ is the $n$-th symmetric power
of the vector space of symmetric matrices.  
In this case $i=0$, so $\nu=\emptyset$. Hence 
$$
\sum\limits_{\mu\vdash n} c_{2\cdot \mu, \emptyset}^\la=
\cases
{1, & \hbox{ if $\la=2\cdot \mu$ for some $\mu\vdash n$ ;}\cr
0, & \hbox{ otherwise. } \cr}
$$
This gives Theorem \ref{Ja}.
Similarly, $\sym^n_n(M_k)= \sym^n(M_k^-)$. In this case $i=n$, $\mu=\emptyset$, 
and a similar computation gives Theorem \ref{Shi}.

\medskip

An analogue of Corollary \ref{dual} follows.

\begin{cor}\label{r.2}
Let $\la,\la'\in\Par_k(2n)$ be conjugate partitions. 
Then, for every $0\le i\le n$,
 the multiplicity
of $V^\la_k$ in $\sym^n_i(M_k)$ is equal to the multiplicity of
$V^{\la'}_k$ in $\sym^n_{n-i}(M_k)$.
\end{cor}

\noindent{\bf Proof.} Combine Theorem \ref{r.1} with identity (2.a).
\qed

\bigskip

\noindent{\bf Proof of Theorem \ref{r.1}.} 
This is a refinement of the proof of Theorem \ref{d.6}.
In this refinement the group $B_n$ appears in an essential way,
whereas in the proof of Theorem \ref{d.6} it was used only as a technical tool.

\smallskip

$M_k^{\otimes n}$ is a $B_n$-module, and $\sym^n(M_k)$ is its submodule, 
for which the $B_n$-action, when restricted to $S_n$, is trivial.
Hence, if the irreducible $B_n$-character $\chi^{\mu,\nu}$ appears in
$\sym^n(M_k)$, then
$$
\langle \chi^{\mu,\nu}\downarrow^{B_n}_{S_n}, 1_{S_n}\rangle \not= 0.
$$
By Lemma \ref{d.2}(a),
$$
\langle \chi^{\mu,\nu}\downarrow^{B_n}_{S_n}, 1_{S_n}\rangle=
\langle \chi^{\mu,\nu}, 1_{S_n}\uparrow^{B_n}_{S_n}\rangle=
\langle \chi^{\mu,\nu}, \sum\limits_{j=0}^n \chi^{(n-j),(j)}\rangle,
$$
and this is nonzero (and equal to 1) if and only if
$\mu=(n-j)$ and $\nu=(j)$ for some $0\le j\le n$.

Combining this with Lemma \ref{r.3} we conclude that $\chi^{(n-i),(i)}$
is the unique irreducible $B_n$-character in $\sym^n_i(M_k)$.

\noindent
Now, as in the proof of Theorem \ref{d.6},
the multiplicity of $V^\la_k$
in $\sym^{n}_i(M_k)$ is 
$$
\langle \chi^\la\downarrow_{B_n}^{S_{2n}}, \chi^{(n-i),(i)} \rangle=
\langle \chi^\la, \chi^{(n-i),(i)}\uparrow_{B_n}^{S_{2n}}\rangle.
$$
By Lemmas \ref{d.3} and \ref{d.0}(a)-(b),
$$
\chi^{(n-i),(i)}\uparrow_{B_n}^{S_{2n}}=
(\chi^{(n-i),\emptyset}\otimes
\chi^{\emptyset,(i)})\uparrow^{S_{2n}}_{B_{n-i}\times B_{i}}=
$$
$$
=(\sum\limits_{\mu\vdash n-i} \chi^{2\cdot \mu} \otimes 
\sum\limits_{\nu\vdash i} \chi^{(2\cdot \nu)'})
\uparrow_{S_{2n-2i}\times S_{2i}}
^{S_{2n}}.
$$
The Littlewood-Richardson rule completes the proof of Theorem \ref{r.1}.

\qed

\subsection{Hook Components of Tensor Powers}

In this subsection we generalize the results of the previous sections
to obtain a bivariate interpolation between symmetric and exterior powers 
of symmetric and skew symmetric matrices. 

\smallskip

As before, the $n$-th tensor power $M_k^{\otimes n}$ carries an $S_n$-action.
The symmetric power $\sym^n(M_k)$ is the $S_n$-invariant part, 
i.e., corresponds to the trivial character $\chi^{(n)}$.
The exterior power corresponds to the sign character $\chi^{(1^n)}$.
We shall denote the factor corresponding to the hook character 
$\chi^{(n-t,1^t)}$ $(0\le t\le n-1)$ by $M_k^{\otimes n}(t)$.
Then

\begin{thm}\label{f.1} 
For every $0\le t\le n-1$ and $\la\in\Par_k(2n)$,
the multiplicity of $V_k^\la$ in $M_k^{\otimes n}(t)$ is
$$
{n-1 \choose t} \sum_{i=0}^{t} (-1)^{t-i} \sigma_{\la}(i) =
{n-1 \choose t} \sum_{i=t+1}^{n} (-1)^{i-t-1} \sigma_{\la}(i)
$$
where
$$
\sigma_{\la}(i) := \sum\limits_{|\alpha|+|\beta|= n-i ,\ |\gamma|+|\delta|=i} 
c_{2\cdot \alpha,(2\cdot \beta)',2*\gamma,(2*\delta)'}^\la 
$$
and the sum runs over 
all partitions $\alpha$ and $\beta$ with total size $n-i$, 
and all partitions $\gamma$ and $\delta$ with distinct parts and total size $i$. 
The operations $\cdot $, $*$ are as defined in Subsection 2.1,
and the extended Littlewood-Richardson coefficients are as defined in
Subsection 2.2.
\end{thm}

\noindent{\bf Proof.}
Similar arguments to those in the proof of Theorem \ref{d.6} show that 
the multiplicity of $V_k^\la$ in the hook component $M_k^{\otimes n}(t)$ is 
$$
{n-1 \choose t} \langle \chi^{(n-t,1^t)} \uparrow^{S_{2n}}_{S_n}, \chi^\la\rangle .
$$
By (3.a),
$$
\langle (\chi^{(n-t,1^t)} + \chi^{(n-t+1,1^{t-1})}) \uparrow^{S_{2n}}_{S_n}, \chi^\la\rangle =
$$
$$
=\langle (\chi^{(n-t)}\otimes \chi^{(1^t)})\uparrow_{S_{n-t}\times S_t}^{S_n}
\uparrow^{S_{2n}}_{S_n}, \chi^\la\rangle =
$$
$$
=\langle (\chi^{(n-t)}\otimes \chi^{(1^t)})
\uparrow_{S_{n-t}\times S_t}^{B_{n-t}\times B_t}
\uparrow_{B_{n-t}\times B_t}^{S_{2n}}, \chi^\la\rangle
$$
By Lemmas \ref{d.2} and \ref{d.3}, for every $t$
$$
(\chi^{(n-t)}\otimes \chi^{(1^t)})
\uparrow_{S_{n-t}\times S_t}^{B_{n-t}\times B_t}=
\chi^{(n-t)}\uparrow_{S_{n-t}}^{B_{n-t}}
\otimes \chi^{(1^t)}
\uparrow_{S_t}^{B_t}=
$$
$$
\sum\limits_{i=0}^{n-t} \chi^{(i),(n-t-i)}\otimes 
\sum\limits_{j=0}^t \chi^{(1^j),(1^{t-j})}=
$$
$$
\sum\limits_{i=0}^{n-t}
(\chi^{(i),\emptyset}\otimes
\chi^{\emptyset,(n-t-i)})
\uparrow_{B_i\times B_{n-t-i}}^{B_{n-t}}
\otimes
\sum\limits_{j=0}^t
(\chi^{(1^j),\emptyset}\otimes
\chi^{\emptyset,(1^{t-j})})
\uparrow_{B_j\times B_{t-j}}^{B_t}.
$$
Hence
$$
(\chi^{(n-t)}\otimes \chi^{(1^t)})\uparrow_{S_{n-t}\times S_t}^{S_{2n}}=
$$
$$
\sum\limits_{i=0}^{n-t}\sum\limits_{j=0}^t
(\chi^{(i),\emptyset} \otimes
\chi^{\emptyset,(n-t-i)} \otimes
\chi^{(1^j),\emptyset} \otimes
\chi^{\emptyset,(1^{t-j})}) 
\uparrow_{B_i\times B_{n-t-i}\times B_{j}\times B_{t-j}}
^{S_{2i}\times S_{2(n-t-i)}\times S_{2j}\times S_{2(t-j)}}
\uparrow_{S_{2i}\times S_{2(n-t-i)}\times S_{2j}\times S_{2(t-j)}}^{S_{2n}}
=
$$
$$
\sum\limits_{i=0}^{n-t}\sum\limits_{j=0}^t
(\chi^{(i),\emptyset} \uparrow_{B_i}^{S_{2i}}
\otimes
\chi^{\emptyset,(n-t-i)} \uparrow_{B_{n-t-i}}^{S_{2(n-t-i)}}
\otimes
\chi^{(1^j),\emptyset} \uparrow_{B_j}^{S_{2j}}
\otimes
\chi^{\emptyset,(1^{t-j})} \uparrow_{B_{t-j}}^{S_{2(t-j)}})
\uparrow_{S_{2i}\times S_{2(n-t-i)}\times S_{2j}\times S_{2(t-j)}}^{S_{2n}}.
$$
Lemma \ref{d.0} and the Littlewood-Richardson rule complete the proof.

\qed

\medskip

\noindent
Let $M_k^{\otimes n}(t,j)$ be the component of $M_k^{\otimes n}(t)$
with $j$ skew symmetric factors.
The following result is a common refinement of Theorems \ref{r.1} and \ref{f.1}.

\begin{thm}\label{hook-sq.1} 
For every $0\le t\le n-1$, $0\le j\le n$ and $\la\in\Par_k(2n)$,
the multiplicity of $V_k^\la$ in $M_k^{\otimes n}(t,j)$ is 
$$
{n-1 \choose t} \sum_{i=0}^{t} (-1)^{t-i} \sigma_{\la}(i,j) =
{n-1 \choose t} \sum_{i=t+1}^{n} (-1)^{i-t-1} \sigma_{\la}(i,j)
$$
where
$$
\sigma_{\la}(i,j) := \sum\limits_{|\alpha|+|\beta|+|\gamma|+|\delta|=n,\ 
|\gamma|+|\delta|=i,\ |\beta|+|\delta|=j}
c_{2\cdot \alpha,(2\cdot \beta)',2*\gamma,(2*\delta)'}^\la ,
$$
and the sum is over all partitions $\alpha,\beta,\gamma,\delta$
with total size $n$ such that $\gamma$ and $\delta$ have distinct parts
and total size $i$, and $\beta$ and $\delta$ have total size $j$.
\end{thm}

\noindent{\bf Proof.}
Lemma \ref{r.3}, used as in the proof of Theorem \ref{r.1}, shows that 
the factors of  $M_k^{\otimes n}(t,j)$, in the decomposition given in 
Theorem \ref{f.1}, are those for which $|\beta|+|\delta|=j$.

\qed

\medskip

\begin{cor}\label{hook-sq.2} 
Let $\la\subseteq (k^k)$ be a partition of $2n$.
For every $0\le t\le n-1$ and $0\le j\le n$,
the multiplicity of $V^\la_k$ in $M_k^{\otimes n}(t,j)$ is equal to 
the multiplicity of $V^{\la'}_k$ in $M_k^{\otimes n}(t,n-j)$.
\end{cor}

\noindent{\bf Proof.} 
By Theorem \ref{hook-sq.1}, it suffices to show that
$$
{n-1 \choose t} \sum_{i=0}^{t} (-1)^{t-i} \sigma_{\la}(i,j) =
{n-1 \choose t} \sum_{i=0}^{t} (-1)^{t-i} \sigma_{\la'}(i,n-j).
$$
This follows from
$$
\sigma_{\la}(i,j) = \sigma_{\la'}(i,n-j),
$$
which in turn follows from (2.a).

\qed

\section{Appendices.}

\subsection{Proof of Lemma \ref{d.2}}

Lemma \ref{d.2} follows from a more general result.

\medskip

\noindent
For partitions $\la=(\la_1,\dots,\la_k)$ and $\mu=(\mu_1,\dots,\mu_m)$,
let $\la\oplus \mu$ be the skew shape defined by 
$$
\la\oplus \mu:= (\la_1+\mu_1, \la_2+\mu_1,\dots,\la_k+\mu_1,\mu_1,\mu_2,\dots,\mu_m)/ (\mu_1^k).
$$

\begin{thm}\label{app.1}
If $(\la,\mu)$ is a bipartition of $n$ then the restriction
$$
\chi^{\la,\mu}\downarrow^{B_n}_{S_n}=\chi^{\la\oplus \mu}.
$$
\end{thm}

\noindent{\bf Proof.}
The characters $\chi^{\la\oplus \mu}$ and
$\chi^{\la,\mu}$, evaluated at elements of $S_n$,
have the same recursive formula (Murnaghan-Nakayama rule).
For $\chi^{\la,\mu}$ see \cite[Theorem 4.3]{Ste2}.
For $\chi^{\la\oplus \mu}$ see \cite[Theorem 5.6.16]{Ke}.

\qed

\noindent{\bf Proof of Lemma \ref{d.2}.}

\noindent
{\bf (a)} Let $(\la,\mu)$ be a bipartition of $n$.
By Frobenius reciprocity and Theorem \ref{app.1},
$$
\langle \chi^{(n)}\uparrow_{S_n}^{B_n}, \chi^{\la,\mu}\rangle=
\langle \chi^{(n)}, \chi^{\la,\mu}\downarrow_{S_n}^{B_n}\rangle=
\langle \chi^{(n)}, \chi^{\la\oplus \mu}\rangle=
\cases
{1, & \hbox{$\max\{\ell(\la),\ell(\mu)\}\le 1$;} \cr
0, & \hbox{otherwise}. \cr}
$$
The last equality follows from the Littlewood-Richardson rule,
reformulated for skew shapes \cite[(7.64)]{Sta}.
By this rule, $\langle \chi^{(n)}, \chi^{\la\oplus \mu}\rangle$
is nonzero (and equal to $1$) if and only if
$\la\oplus \mu$ is a horizontal strip (i.e., each column contains at most one box).

\smallskip

\noindent{\bf (b)} The proof for $\chi^{(1^n)}$ is similar.

\qed

\subsection{Proof of Lemma \ref{r.3}}

\noindent{\bf Proof.}
Let $\sigma_i:=(i,-i)\in B_n$ $(1\le i\le n$),
and let $\eta$ be the sum
$\sum\limits_{i=1}^n \sigma_i\in \bbc[B_n]$.
Consider the tensor product $w=w_1\otimes w_2\otimes \cdots \otimes w_n \in M_k^{\otimes n}$, where each $w_i$ is either a symmetric or a skew symmetric matrix. Then according to the $B_n$-action, defined in Section 5.2,
$$
\sigma_i (w) =\cases
{w, & \hbox{ if $w_i$ is symmetric;} \cr
-w, & \hbox{ if $w_i$ is skew symmetric.} \cr}
$$
Hence, for every vector $v\in M_k^{\otimes n}(i)$
$$
\eta (v) = (n-2i) v . \leqno{(6.a)}
$$

On the other hand,
the set
$\{\sigma_i\ |\ 1\le i\le n\}$ is a conjugacy class in $B_n$. 
Thus the element $\eta=\sum\limits_{i=1}^n \sigma_i$ is in the center of 
$\bbc[B_n]$.
By Schur's Lemma, for every vector $v$ in the irreducible $B_n$-module $S^{\mu,\nu}$
$$
\eta (v)= c^{\mu,\nu}\cdot v, 
$$
where 
$$
c^{\mu,\nu}= {\chi^{\mu,\nu}(\eta)\over \chi^{\mu,\nu}(id)}=
 {n \chi^{\mu,\nu}(\sigma_1)\over \chi^{\mu,\nu}(id)}.
$$

Let $f^\la$, $f^{\mu,\nu}$ be the number of standard Young tableaux 
(bitableaux) of shapes
$\la$, $(\mu,\nu)$ respectively.
Recall that 
$$
\chi^{\mu,\nu}(id)=f^{\mu,\nu}={n\choose |\nu|}f^\mu f^\nu,
$$
and that 
$\chi^{\mu,\nu}(\sigma_1)$ is equal to the number of standard Young bitableaux of shape $(\mu,\nu)$, in which the digit 1 is in the first diagram $\mu$, minus the number of those in which 1 is in
the second diagram $\nu$.
Thus
$$ 
\chi^{\mu,\nu}(\sigma_1)={n-1\choose |\nu|}f^\mu f^\nu-{n-1\choose |\nu|-1}f^\mu f^\nu=
{n-2|\nu|\over n}{n\choose |\nu|}f^\mu f^\nu.
$$
It follows that
$$
c^{\mu,\nu}=
{n \chi^{\mu,\nu}(\sigma_1)\over \chi^{\mu,\nu}(id)}=n-2|\nu|,
$$
and therefore
$$
\eta(v)= (n-2|\nu|)v \qquad (\forall v \in S^{\mu,\nu}). \leqno{(6.b)}
$$
Combining (6.a) with (6.b) completes the proof.

\qed

\bigskip

\noindent{\bf Acknowledgments.} The authors thank R.\ Howe, A.\ Regev, T.\ Seeman, N.\ Wallach
and an anonymous referee for their useful comments.


\begin{thebibliography}{xx}

\bibitem{Dv} Y. Dvir, {\it On the Kronecker product of $S_n$
characters.} J. Algebra 154 (1993), 125--140.

\bibitem{Eh} C. Ehresmann, {\it Sur la topologie de certains espaces
homog\`enes.} Ann. of Math. 35 (1934), 396--443.

\bibitem{GW} R. Goodman and N. R. Wallach,
{\it Representations and Invariants of the Classical Groups.}
Encyclopedia of Math. and its Appl. Vol. 68,
Cambridge University Press, 1998.

\bibitem{He}S. Helgason, {\it A duality for  symmetric spaces with applications.} Adv. Math. 5 (1970), 1--54.

\bibitem{H} R. Howe,
{\it Remarks on Classical
Invariant Theory}.
Trans. Amer. Math. Soc. 313 (1989), 539--570.

\bibitem{Ho} R. Howe and T. Umeda, {\it The Capelli identity, the double commutant theorem, and multiplicity-free actions.} Math. Ann. 290 (1991), 565--619.

\bibitem{Ja} A. T. James, {\it Zonal polynomials of the real positive definite matrices.} Annals of Math. 74 (1961), 456--469.

\bibitem{JK} G. D. James and A. Kerber, {\it The Representation Theory of the Symmetric Group.} Encyclopedia of Math. and its Appl. Vol. 16,
      Addison-Wesley, 1981.

\bibitem{Ke} A. Kerber, {\it Algebraic combinatorics via finite group actions.} Bibliographisches Institut, Mannheim, 1991. 

\bibitem{Ko} B. Kostant, {\it Lie algebra cohomology and the generalized Borel Weil theorem.} Ann. of Math. 74 (1961), 329--387.

\bibitem{Se} J. P. Serre, {\it Linear Representations of Finite Groups.} 
 Springer-Verlag, 1977.

\bibitem{Md} I. G. Macdonald, 
{\it Symmetric Functions and Hall Polynomials}. second edition, Oxford Math.\ Monographs, 
Oxford Univ.\ Press, Oxford, 1995.

\bibitem{Re} A. Regev, {\it The Kronecker product of $S_n$
characters and an $A\otimes B$ theorem for Capelli identities.} J. Algebra 66 (1980), 505--510.


\bibitem{Sh} G. Shimura, {\it On differential operators attached to certain representations of classical groups.} Invent. Math. 77 (1984), 463--488.

\bibitem{Sta} R. P. Stanley, {\it Enumerative Combinatorics, Volume II.}
Cambridge Univ. Press, Cambridge, 1999.

\bibitem{St} J. R. Stembridge, {\it On Schur $Q$-functions and the primitive idempotents of commutative Hecke algebra.} J. Alg. Combin.
1 (1992), 71--95.

\bibitem{Ste2} 
J.\ Stembridge, 
{\it On the eigenvalues of representations of reflection groups and wreath 
products.} 
Pacific J.\ Math.~140 (1989), 359--396.

\bibitem{Th} R. Thrall, {\it On symmetrized Kronecker powers and the structure of the free Lie ring.} Amer. J. Math. 64 (1942), 371--388.

\end{thebibliography}
\end{document}